%Paper: q-alg/9502004
%From: gannon@abacus.concordia.ca (Gannon Terry)
%Date: Thu, 2 Feb 95 19:38:53 -0500

\magnification=\magstep1  \overfullrule=0pt
\def\la{\lambda} \def\om{\omega} \def\sig{\sigma} \def\eps{\epsilon}
\def\al{\alpha}  \def\Eqde{\,\,{\buildrel {\rm def} \over =}\,\,}
\def\be{\beta}    \def\r{\bar{r}}    \def\A{A^{(1)}} \def\La{\Lambda}
\font\huge=cmr10 scaled \magstep2 \def\Rh{\overline{\rho}}
\def\h{h^\vee} \def\ga{\gamma}

\rightline{September, 1994}
\bigskip\bigskip
\centerline{{\bf \huge Symmetries of the Kac-Peterson}}
\bigskip \centerline{{\bf \huge Modular Matrices of Affine Algebras}}\bigskip
\bigskip\bigskip
\centerline{Terry Gannon\footnote{$^1$}{Present address: Math Dept, Concordia
University, Montr\'eal, Canada, H3G 1M8}}
\bigskip\centerline{{\it Institut des Hautes Etudes Scientifiques}}
\centerline{{\it 91440 Bures-sur-Yvette, France}}\bigskip\bigskip
\bigskip\bigskip
\noindent{{\bf 0\qquad Introduction}} \medskip

One of the more fruitful observations in recent years was the discovery [19]
that the (normalized)
character $\chi_\mu$ of the integrable highest weight $g$-module
$L(\mu)$ of a nontwisted Kac-Moody algebra $g=X^{(1)}_r$ is a modular form.
Thanks to the structure of the affine Weyl group, these $\chi_\mu$ can be
written as ratios of theta functions, so on them exists a natural action
of $SL_2({\bf Z})$. In fact, this action defines a representation
$$R:SL_2({\bf Z})\rightarrow GL_{n}({\bf C})$$
on the space spanned by the $\chi_\mu$, for $\mu$ lying in the set $P_+(g,k)$
of all level $k$ highest weights of $g$ (the dimension $n$ equals the
cardinality $\|P_+(g,k)\ {\rm mod}\ {\bf C}\delta\|$, for the imaginary root
$\delta$ of $g$). The matrices in the image
$R(SL_2({\bf Z}))$, which we will call the {\it Kac-Peterson modular
matrices}, are unitary (in the $\chi_\mu$ basis).

This observation and its various consequences have
found application in several areas. For example, many conformal field
theories [29,14,15] are intimately connected with nontwisted affine
$g$.
Much of their structure is encoded in their ``genus 1 partition
function'', which in the case of a Wess-Zumino-Witten theory [29,14] at
level $k$ based on $g$ is a sesquilinear combination of
characters of $g$:
$$Z=\sum_{\mu,\nu\in P_+(g,k)\ {\rm mod}\ {\bf C}\delta}M_{\mu,\nu}\,
\chi_\mu\,\chi_\nu^*.\eqno(1)$$
This function must satisfy some properties, namely:

\item{$\bullet$} the coefficients $M_{\mu,\nu}$ in eq.(1) must be non-negative
integers;

\item{$\bullet$} $Z$ must be invariant under the action of $SL_2({\bf Z})$ on
the characters -- equivalently, the matrix $M$ must commute with all
matrices in $R(SL_2({\bf Z}))$;

\item{$\bullet$} $M_{k\Lambda_0,k\La_0}=1$.

We shall call any such matrix $M$ a {\it physical invariant}. The
classification
of these physical invariants is a major problem in conformal field theory.
Unfortunately, it is difficult and in spite of much effort little is
known. What makes the problem more interesting though is that many startling
coincidences have appeared between the few existing classifications
and other areas of mathematics.
For instance, the physical invariants corresponding to $g=A_1^{(1)}$
[5] fall into the A-D-E pattern familiar to e.g.\ singularity theory, the
finite subgroups of $SU_2({\bf C})$, and of course the simply-laced Lie
algebras. In the classification of the physical invariants for $g=
A_2^{(1)}$ [10] a number of surprising relations [23,7] with the Jacobians of
the Fermat curves [22] have been found. Recently [30], a more inclusive and
sophisticated interpretation of some of these relationships has been proposed
using generalized Coxeter graphs.

It should not be too surprising that as rich a subject as conformal
field theory has subtle interconnections with other parts of mathematics. And
indeed some have become established in recent years, thanks to the
work of Witten, Kac, Verlinde, and many others.
But the coincidences involving the CFT classifications are still poorly
understood; it is difficult to know their significance or to anticipate
how they extend to other affine algebras. For this reason, as
well of course for the classification of conformal field theories itself,
it is certainly of interest to try to find all physical invariants
corresponding to e.g.\ $\A_r$.

In Section 2 we classify an important subset of the $\A_r$ physical
invariants, called {\it automorphism invariants} (see eqs.(5)). These
correspond to those $M$ in eq.(1) which are permutation matrices;
equivalently, they are the exact symmetries of all the Kac-Peterson modular
matrices of the affine algebra at a given level. We extend our proof to
$\A_{\vec{r}}\Eqde(A_{r_1}\oplus\cdots\oplus A_{r_s})^{(1)}$, $\vec{r}=
(r_1,\ldots,r_s)$ in Section 3; it should also
extend to any other affine algebra, as
well as a much larger class of physical invariants (namely all automorphisms
of the Bernard [2] -- i.e.\ simple current [25] -- chiral extensions).
The value of our result to the classification problem of CFT will be
discussed in more depth at the end of Section 1.

A brief sketch of the proof is given at the beginning of Section 2.
Among the tools we use in this paper are the Verlinde formula, and the fact
that the Weyl character of any representation of a finite-dimensional complex
simple Lie algebra $\bar{g}$
can be written as a polynomial in the Weyl characters of its fundamental
representations. An important connection between the representation theory
of $g$ and that of finite dimensional $\bar{g}$ is the fact (see eq.(3a))
that certain ratios of entries of the Kac-Peterson modular matrix $S$ equal
certain values of the Weyl characters of $\bar{g}$. This connection is
exploited throughout the paper.

In the remainder of this introduction we will briefly review the little
that is known about physical invariants and their classification.

It is easy to show [9] that for any fixed algebra $g$ and level $k$,
there  will be only finitely many physical invariants. Many of these have
already been constructed. Almost all known ones are built up from the
symmetries $Aut\ \Pi^\vee$ of the
extended Dynkin diagram, in simple ways [2,1,25]; they are called
${\cal D}$-{\it type}
invariants (by analogy with the A-D-E classification of $\A_1$).
The remaining physical invariants are called {\it exceptional} (see e.g.\
[4,27]). Most known exceptionals are related to the so-called {\it conformal
subalgebras} of $g$ [4,20] -- i.e.\ subalgebras
whose branching functions $b^\la_\mu(\tau)$ are constants.

The best known classification
is the A-D-E one for $\A_{1}$ [5]. There we find an ${\cal A}_k$-series
(for all $k$) and ${\cal D}_k$-series (for even $k$) -- both these
are of ${\cal D}$-type. There are also three exceptionals, at levels 10, 16
and 28. The only other classifications known at present for the physical
invariants of affine algebras $g$ are: $k=1$ when $\bar{g}$ is simple [9,16];
$A_2^{(1)}$ for all $k$ [10]; and $(A_1\oplus A_1)^{(1)}$ for all $k=(k_1,
k_2)$ [11].

We will prove that the only automorphism invariants for $\A_{\vec{r}}$ are
those of ${\cal D}$-type. Incidently, this is not the case for all
affine algebras: e.g.\ there are exceptional automorphism invariants for
$G_{2}^{(1)}$ at level 4 and $F_{4}^{(1)}$ at level 3 [27], and for
infinitely many $B_r^{(1)}$ and $D_r^{(1)}$ at level 2 [12]. The methods
developed here should apply to the other affine algebras [12], though the
details will certainly be messier.

The only other existing major classification of automorphism invariants
is for the ``simple current automorphism invariants'' [13]
of any rational conformal field theory
(subject to a condition on the corresponding modular matrix $S$).
However the argument in [13] cannot apply
here because it limits its attention to automorphism invariants of a special
form not shared by most of the $\A_{\vec{r}}$
automorphism invariants; it also assumes that (roughly speaking)
the modular matrix $S$ in eq.(2a) does not have too many zeros --
something known at present only for $A_1^{(1)}$ and $A_2^{(1)}$.

Years of effort with little accomplished has made most in the field
skeptical about the possibility of a classification of all physical
invariants, even for $A_r^{(1)}$. Our result, especially considering its
simplicity, directly challenges this pessimism. Knowing all the
automorphism invariants is a necessary and major
step towards the full classification.

\bigskip \noindent{{\bf 1 \quad The Kac-Peterson modular matrices}}\medskip

In this section we will introduce some notation and terminology, and review
some results in [19,20]. Our notation will remain as close as possible to
[17]. We will focus here on $g=\A_r$; analogous statements hold for the other
affine algebras [17,21], and those for $\A_{\vec{r}}$ will be given at the
beginning of Section 3.

Write $\r=r+1$ and $\bar{k}=k+\r$. Choose a Cartan subalgebra $h$ for
$g$; its dual $h^*$ will be spanned by the fundamental weights
$\La_0,\ldots,\La_r$ together with the imaginary root $\delta=\sum_{i=0}^r
\alpha_i$. We can (p.155 of [17], eq.(1.5.12) of [20]) and will identify a
highest weight
$$\la=\sum_{i=0}^r\la_i\La_i+z\delta\in h^*$$
with its Dynkin labels $(\la_0,\la_1,\ldots,\la_r)$, and drop the $z$.
In particular, the set of level $k$ highest weights for $\A_r$ becomes
$$P_+^{r,k}\Eqde P_+(\A_r,k)\ {\rm mod}\ {\bf C}\delta
=\bigl\{\sum_{i=0}^r\la_i\La_i\,|\,\la_i\in{\bf Z},\
\la_i\ge 0,\ \sum_{i=0}^r\la_i=k\,\bigr\}.$$
Put $\rho=\sum_{i=0}^r\La_i$. The invariant bilinear form $(-|-)$ for
the underlying finite-dimensional algebra $\bar{g}=A_r$,  normalized
so that the roots have norm 2, can be extended to $h^*$ by defining
$(\La_0|\La_i)=(\La_i|\La_0)=(\delta|\delta)=0$, $(\delta|\Lambda_i)=(\Lambda_i
|\delta)=\delta_{i,0}$, $\forall i$. We will  let $\overline{\la}$
denote the orthogonal projection $(\la_1,\ldots,\la_r)$ onto the dual
$\bar{h}^*={\bf C}\Lambda_1\oplus\cdots\oplus{\bf C}\Lambda_r$
of the Cartan subalgebra  of   $\bar{g}$.
 Let $Aut\,\Pi^\vee$ denote the group of automorphisms of the
(extended) Dynkin diagram of $g$: for $\A_r$ it has generators $C_r$
and $J_r$ acting on $P_+^{r,k}$ by
$$\eqalign{C_r\la&=\la_0\,\La_0+\sum_{i=1}^r\la_{\r-i}\,\La_i,\cr
J_r\la&=\la_r\,\La_0+\sum_{i=1}^r\la_{i-1}\,\La_i.\cr}$$

Consider the irreducible integrable highest weight $g$-module $L(\la)$.
 Let $L(\la)=\oplus_\beta
L(\la)_\beta$ denote its weight space decomposition with respect to
$h$. Define the normalized character of $L(\la)$ to be
$$\chi_\la=e(-m_\la\delta)\,\sum_{\beta}dim\, \bigl[L(\la)_\beta\bigr]\,\,
e(\beta),$$
where $m_\la$ is a rational number (the {\it modular anomaly}), and
for each $\gamma\in h^*$, $e(\gamma)$ can be taken to be the function
$e(\gamma):h\rightarrow {\bf C}$ given by $e(\gamma)(v)=e^{\gamma(v)}$.

Coordinatizing $h$ appropriately [19], $\chi_\la$ may be regarded as a
complex-valued function of $\bar{z}\in \bar{h}$ and complex variables $\tau,
u$, such that for each $\la\in P_+^{r,k}$
$$\eqalignno{\chi_\la({-1\over\tau},{\bar{z}\over\tau},u-{(\bar{z}|\bar{z})
\over 2\tau})=&\sum_{\mu\in P_+^{r,k}}S^{r,k}_{\la,\mu}\,\,\chi_\mu,&\cr
S^{r,k}_{\la,\mu}\Eqde\,&{\exp[\pi i r\r/4]\over \bar{k}^{r/2}\,\sqrt{\r}}\,
\sum_{w\in \overline{W}}det(w)\,\exp[-2\pi i{\bigl(w(\overline{\la}+\Rh)\,|\,
\overline{\mu}+\Rh\bigr)\over \bar{k}}];&(2a)\cr
\chi_\la(\tau+1,\bar{z},u)=&\sum_{\mu\in P_+^{r,k}}T^{r,k}_{\la,\mu}\,\,
\chi_\mu, &\cr
T^{r,k}_{\la,\mu}\Eqde\,&\exp[\pi i\{{(\la+\rho\,|\,\la+\rho)\over
\bar{k}}-{(\rho\,|\,\rho)\over \h}\}]\,\,\delta_{\la,\mu}.&(2b)\cr}$$
We will often delete the superscripts `$r,k$' in the following.
$\overline{W}$ in eq.(2a) denotes the (finite) Weyl group of $A_r$.
The transformations in eqs.(2) generate the complete action of $SL_2({\bf
Z})$ on the space spanned by the $\chi_\la$,
and so these equations suffice to uniquely specify the representation
$R$ of $SL_2({\bf Z})$. The Kac-Peterson matrices $R(SL_2({\bf Z}))$ consist
of all possible products of $S^{r,k}$ and $T^{r,k}$.

The matrices $S^{r,k}$ and $T^{r,k}$ are unitary and symmetric. $(S^{r,k})^2
=C^{r,k}$, the matrix characterizing the action of $C_r$ on $P_+^{r,k}$. A
remarkable property [19] of $S^{r,k}$, which we will use frequently, follows
immediately from eq.(2a) and the Weyl character formula:
$${S_{\la,\mu}\over S_{k\La_0,\mu}}=ch_{\overline{\la}}
\bigl(-2\pi i \,\nu^{-1}(\overline{\mu}+\Rh)\,/\bar{k}\bigr),\eqno(3a)$$
where $\nu:\bar{h}\rightarrow \bar{h}^*$ is the isomorphism defined by $(-|-)$,
and where $ch_{\overline{\la}}$ denotes the Weyl character of the $A_r$-module
$\overline{L}(\overline{\la})$:
$$ch_{\overline{\la}}\Eqde\sum_{\overline{\beta}}dim\,\bigl[\overline{L}(
\overline{\la})_{\overline{\beta}}\bigr]\,\,e(\overline{\beta}).$$
A special case of eq.(3a) is the {\it q-dimension} of any $\la\in P^{r,k}_+$:
$$Q^{r,k}(\la)\Eqde{S_{\la,k\La_0}\over S_{k\La_0,k\La_0}}=
\prod_{\overline{\alpha}>0}
{\sin[\pi\,(\overline{\la}+\Rh\,|\,\overline{\alpha})/\bar{k}]\over \sin[
\pi\,(\overline{\rho}\,|\,\overline{\al})/\bar{k}]}.\eqno(3b)$$
The $\overline{\alpha}>0$ in eq.(3b) are the positive roots of $A_r$. These
considerations also imply [17]
$$S_{k\La_0,\la}\ge S_{k\La_0,k\La_0}>0,\quad \forall\la\in P_+^{r,k}.\eqno(3c)
$$

Straightforward calculations give us [20]
$$\eqalignno{t(\la)\Eqde &\,\r\,\,(\la\,|\,\La_r)=\sum_{j=1}^r j\la_j,&(4a)\cr
t(J_r^a\la)\equiv &\, ka+t(\la)\quad ({\rm mod}\ \r),&(4b)\cr
T_{J_r^a\la,J_r^a\mu}^{r,k}=& \,\exp[\pi i\,\bigl(-2a\,t(
\la)+ka\,(\r-a)\bigr)/\r]\, \,T^{r,k}_{\la,\mu},&(4c)\cr
S_{J_r^a\la,J_r^b\mu}^{r,k}=&\,\exp[2\pi i\,\bigl(b\,t(\la)+a\,t(\mu)+kab\bigr)
/\r]\,\,S_{\la,\mu}^{r,k},&(4d)\cr
S_{C_r\la,\mu}^{r,k}=&\,S^{r,k}_{\la,C_r\mu}=S^{r,k*}_{\la,\mu},\qquad
T^{r,k}_{C_r\la,C_r\mu}=T^{r,k}_{\la,\mu}.&(4e)\cr}$$

Incidently, there is a Galois symmetry [9,23,6] obeyed by $S^{r,k}$
which has proven to be valuable in some contexts
(see e.g.\ [10]). Although we will not use it here,
it is sufficiently little known to warrant us repeating.
We see from eq.(2a) that the entries of $S^{r,k}$ lie in the cyclotomic
extension ${\bf Q}(\exp[\pi i/2\bar{k}\r])$ of the rationals ${\bf Q}$.
Choose any element $g$ of the corresponding
Galois group. Then associated to $g$ is a permutation $\la\mapsto \la^g$
of the weights in $P_+^{r,k}$, as well as a map $\eps_g:P_+^{r,k}
\rightarrow\{\pm 1\}$, such that
$$g(S_{\la,\mu})=\eps_g(\la)\,S_{\la^g,\mu}=\eps_g(\mu)\,S_{\la,\mu^g},
\qquad \forall \la,\mu\in P_+^{r,k}.$$
The identical result holds for any other nontwisted
affine algebra. This Galois symmetry reappears throughout rational
conformal field theory: for instance the matrix $M$ in eq.(1) must obey
$$M_{\la,\mu}=\eps_g(\la)\,\eps_g(\mu)\,M_{\la^g,\mu^g},\qquad \forall
\la,\mu\in P_+(g,k)\ {\rm mod}\ {\bf C}\delta.$$

\noindent{\bf Definition.}\quad {\it An {\bf automorphism invariant} of $g$
at level $k$ is a permutation $\sigma$ of $P_+(g,k)$ mod ${\bf C}\delta$
such that $U_{\la,\mu}=
U_{\sigma\la,\sigma\mu}$ for all Kac-Peterson matrices $U\in R(SL_2({\bf Z}))$
and all $\la,\mu\in P_+(g,k)$ mod ${\bf C}\delta$.} \medskip

\noindent{{\it Remarks.}}\quad
Since any such $U$ is generated by the analogues of the matrices $S$ and $T$
in eqs.(2), it suffices to require that $\sigma$ satisfy
$$\eqalignno{S_{\la,\mu}=&S_{\sig\la,\sig\mu},\qquad \forall \la,
\mu\in P_+(g,k)\ {\rm mod}\ {\bf C}\delta,&(5a)\cr
T_{\la,\mu}=&T_{\sig \la,\sig\mu},\qquad \forall
\la,\mu\in P_{+}(g,k)\ {\rm mod}\ {\bf C}\delta.&(5b)\cr}$$

Our task in this paper is to find all automorphism invariants for
$\A_{\vec{r}}$. The simplest examples are
$\sig=C_{r_1}^{c_1}\otimes\cdots \otimes C_{r_s}^{c_s}$ for any choice of
$c_i\in {\bf Z}$ -- these obey eqs.(5) because of eq.(4e).
Of course the set of all automorphism invariants will be a group under
composition. Because $k\La_0$ is the only row of $S$ which is
strictly positive, we see that
$$\sig(k\La_0)=k\La_0.$$

Defining $M$ by $M_{\la,\mu}=\delta_{\mu,\sig\la}$, we find that {\it any
automorphism invariant is a physical invariant}. The converse is usually not
true; nevertheless the automorphism invariants are an important subset
of the physical invariants which until now were quite intractible.
For one thing, the matrix product of an automorphism invariant with any
other physical invariant will also be a physical invariant. Also, knowing
the automorphism invariants for affine algebras automatically means we
know all automorphisms of those exceptional chiral extensions due to
conformal subalgebras [4,20] -- the main source of exceptional extensions.
Finally, knowing the automorphism invariants should permit the
classification of all physical invariants $M$ with the property:
$$M_{\la,k\La_0},\,M_{k\La_0,\mu}\ne 0\ \Rightarrow\ \exists J,J'\in Aut\,
\Pi^\vee\ {\rm such\ that}\ \la=J(k\La_0),\,\mu=J'(k\La_0).$$
This would be important because almost every known physical invariant
obeys  that property. Of course, any automorphism invariant satisfies it
-- in fact a physical invariant $M$ will be an automorphism
invariant iff $M_{\la,k\La_0}=M_{k\La_0,\la}=\delta_{\la,k\La_0}$ $\forall
\la\in P_+(g,k)$ mod ${\bf C}\delta$ [9].

\bigskip\noindent{{\bf 2 \quad The automorphism invariants of $\A_r$}}\medskip

Our goal in this paper is to find the automorphism invariants
of $\A_{\vec{r}}$, $\vec{r}=(r_1,\ldots,r_s)$. The main ideas however are
present in the much simpler special case $s=1$. In order to make the general
argument clearer, in this section we will limit the discussion to $\A_r$.

The proof is surprisingly simple. Calculus will show that the
q-dimensions eq.(3b) force the automorphism invariant $\sig$ to map
the weight $\om_1=(k-1)\La_0+\La_1$ to $C^a_rJ^b_r\om_1$, for some $a,b$.
Using eq.(5b) we will see that $\sig\om_1=C_r^a\sig_m\om_1$ for some
known automorphism invariant $\sig_m$. Thus it suffices to consider the
case where $\sig\om_1=\om_1$. Using Verlinde's formula, we can then
show that $\sig$ must fix all $\om_i=(k-1)\La_0+\La_i$. But any Weyl
character $ch_{\overline{\la}}$ is a polynomial in the $ch_{\overline{\om_i}}$
-- from eq.(3a) this tells us $\sig$ must fix all weights.

Because of eqs.(4c),(4d) we can expect to build automorphism invariants from
$J_r$. In particular [25], define
$$\tilde{k}=\left\{\matrix{\bar{k}&{\rm if}\ r\equiv \bar{k}\equiv 0\quad {\rm
(mod\ 2)};
\cr k&{\rm otherwise}.\cr}\right. \eqno(6a)$$
Choose any positive integer $m$ dividing $\r$ such that
$m\tilde{k}$ is even and $gcd\{\r/m,m\tilde{k}/2\}=1$. Then  we
can find an integer $v$ such that $vm\tilde{k}/2\equiv 1$ (mod $\r/m$).
To each such divisor $m$ of $\r$, we can define an automorphism invariant
$\sig_m$ given by
$$\sig_m\la\Eqde J_r^{-v\,m\,t(\la)}\la.\eqno(6b)$$
For example, $\sig_{\r}$ corresponds to the identity permutation.
That each $\sig_m$ is a bijection of $P_+^{r,k}$, follows from $\sig_m\circ
\sig_m=id.$. That they satisfy eqs.(5) can be verified explicitly, using
eqs.(4).

\medskip\noindent{\bf Theorem 1}\quad {\it The only automorphism invariants
$\sig$ of $\A_r$ at level $k$ are $C_r^a\sig_m$, for $a=0,1$ and $\sig_m$
defined in eq.(6b). All of these are distinct, except for $r=1$ (when
$C_r=id.$) or $k\le 2$.}\medskip

\noindent{{\it Remark 1}} \quad From Theorem 1 we find that there are
precisely $2^{c+p+t}$ automorphism invariants for $\A_{r}$ at level $k$, where:

\item{$c$}$=1$, unless $r=1$ or $k\le 2$ when $c=0$ (an exception: $c=-1$ when
both $r=1$ and $k=2$);

\item{$p$} is the number of distinct odd primes which divide $\r$ but not $k$;

\item{$t$}$=0$ if either $r$ is even, or $r$ is odd and $k\equiv 0$
(mod 4), or $k$ is odd and $r\equiv 1$ (mod 4) -- otherwise $t=1$.

\medskip\noindent{{\it Remark 2}}\quad
{}From the theorem we also see that all automorphism invariants are of order
2 (i.e.\ $\sig^2=id.$), and commute: in fact
$$C_r^a\sig_m\circ C_r^b\sig_{m'}=C_r^{a+b}\sig_{m''},\ {\rm where}\ m''=\r \,
gcd\{m,m'\}^2/mm'.$$
Both of these facts are surprising, and not true for general $g$ (as
we will see in Section 3).

\medskip Three special cases of the theorem were known
previously: $r=1$ [5], $k=1$ [16], and $r=2$ [10,24]. The remainder of this
section is devoted to the proof of Theorem 1.

For each $i=1,\ldots,r$ define $\om_i=(k-1)\La_0+\La_i$. For any $\la$
define the $Aut\,\Pi^\vee$-orbit $[\la]=\{C_r^aJ_r^b\la\,|\,a,b
\in{\bf Z}\}$. Let $o(\la)$
equal the number of indices $0\le i\le r$ such that $\la_i>0$. For example,
$\la\in[k\La_0]$ iff $o(\la)=1$. Also, $o(\om_i)=2$ (provided $k>1$).

Let us begin by proving the second statement in the Theorem.
That the $\sig_m$'s are all distinct is easy to see by looking at $\sig_m
\om_1$. Next, for $r>1$ and $k>2$, $C_r\om_1\not\in\{J_r^i\om_1\,|\,i\in{\bf
Z}\}$, so $C_r\sig_m
\om_1\ne\sig_{m'}\om_1$ for any $m,m'$. Incidently, when $k\le 2$,
$C_r=\sig_m$ where $m=1$ or 2.

Recall $Q^{r,k}(\la)$ defined in eq.(3b). It will be constant within each
orbit $[\la]$. Eq.(3b) was analysed in [8], by extending the
domain  of $Q$ to {\it real} (as opposed to {\it integer}) vectors $\la$. We
will make use of this idea below.

\medskip \noindent{{\bf Proposition 1}}\quad {\it For $k>1$ and any $\la
\not\in[k\La_0]\cup[\om_1]$, $Q^{r,k}(\la)>Q^{r,k}(\om_1)>Q^{r,k}(k\La_0)$.}

\noindent{{\it Proof.}}\quad
Consider any $\la\in P_{+}^{r,k}$, $\la\not\in[k\La_0]$. Suppose first
that $o(\la)\ge 3$, and let $\la_i,\la_j>0$ for $i\ne j>0$. Consider $\mu(t)=
\la+t\La_j-t\La_i$. Then $(\mu(\la_i))_i=(\mu(-\la_j))_j=0$. An easy
calculation, similar to one in [8], tells us
$${{\rm d}\over {\rm d}t}Q(\mu(t))\bigl|_{t=t_0}=0\,\Rightarrow\,{{\rm d}^2
\over
{\rm d}t^2}Q(\mu(t))\bigr|_{t=t_0}=-Q(\mu(t_0))\,{\pi^2\over \bar{k}^2}\,
\sum_{\overline{\alpha}>0}{(\overline{a}|\overline{
\alpha})^2\over \sin^2[\pi (\overline{\mu(t_0)+\rho}| \overline{\alpha})/
\bar{k}]}<0, $$
where $\overline{a}=\overline{\La_j}-\overline{\La_i}$. Thus
$Q(\mu(t))$ will take its minimum on the endpoints, i.e.\ either
$Q(\mu(\la_i))<Q(\la)$ or $Q(\mu(-\la_j))<Q(\la)$. In either case,
we have found a $\mu\in P_{+}^{r,k}$ with $Q(\mu)<Q(\la)$, which has
$o(\mu)=o(\la)-1$.

Continuing inductively, we see it suffices to consider the weights with
$o(\la)=2$. The same argument allows us to put one of those two Dynkin
labels equal to 1. In other words, starting with any weight $\la\not\in
[k\La_0]$, we can find a $\om_\ell$ such that $Q(\la)\ge Q(\om_\ell)$, with
equality iff $\la\in[\om_\ell]$. All that remains is to compare $Q(\om_\ell)$
with $Q(\om_1)$.

We find from eq.(3b) that for all $1<m\le r$,
$${Q(\om_{m})\over Q(\om_{m-1})}=\prod_{i=1}^{m-1}{\sin[\pi{m-i
\over\overline{k}}]\over\sin[\pi{m+1-i\over \overline{k}}]}\,\prod_{j=m}^{r}{
\sin[\pi{j+2-m\over\overline{k}}]\over\sin[\pi{j+1-m\over \overline{k}}]}={
\sin[\pi{r+2-m\over\overline{k}}]\over\sin[\pi{m\over\overline{k}}]}. $$
Hence $Q(\om_1)<Q(\om_\ell)$ unless $\ell=1$ or $r$. But $\om_{r}=C_r\om_1\in
[\om_1]$. Finally, compute $Q(\om_1)=\sin[\pi \r/\bar{k}]/\sin[\pi/\bar{k}]
\ge 1= Q(k\La_0)$, with equality only at $k=1$.\qquad  QED \medskip

The restriction in Prop.\ 1 to $k>1$ is not important, because $P_+^{r,1}=
[k\La_0]$. Prop.\ 1 together with eq.(5a) tells us, for any $r\ge 1$, $k\ge
1$, that $\sig\om_1=C_r^aJ_r^b\om_1$ for some $a,b$. By replacing $\sig$ with
$C_r^a\sig$, we may assume $a=0$. Evaluating the expression $T^{r,k}_{\om_1,
\om_1}=T^{r,k}_{J_r^b\om_1,J_r^b\om_1}$ using eq.(4c), we get
$$2b\equiv kb\,(\r-b)\quad ({\rm mod}\ 2\r).\eqno(7a)$$

\medskip\noindent{{\bf Proposition 2}}\quad {\it Let $b$ be any integer
satisfying
eq.(7a). Then $\sig_m(J_r^b\om_1)=\om_1$, for one of the $\sig_m$ of eq.(6b).}

\noindent{{\it Proof.}}\quad Given any solution $b$ to eq.(7a), define $m=gcd\{
b,\r\}$ and $v'=b/m$ -- we can suppose (adding
$\r$ to $b$ if necessary) that $v'$ will be coprime to $2\r/m$.
Suppose we can show
$$\tilde{k}b/2 \equiv -1\quad ({\rm mod}\ \r/m).\eqno(7b)$$
Then $\tilde{k}m/2$ will be an integer coprime
to $\r/m$ (so that $\sig_m$ exists), and $b\equiv -vm$ (mod $\r$) for $v$
as in eq.(6b) (so that $\sig_m(J_r^b\om_1)=\om_1$).

Now, dividing eq.(7a) by $b$ gives
$$kb\equiv -2+k\bar{r}\quad ({\rm mod}\ 2\r/m).\eqno(7c)$$
For $\r$ odd, $b$ is odd and eq.(7b) is immediate from eq.(7c). For $\r$ even,
eq.(7c) requires either $k$ or $m$ to be even, so again eq.(7b) follows from
eq.(7c). \qquad QED  \medskip

By Prop.\ 2 we see that, replacing $\sig$ with $\sig_m\circ\sig$, {\it it
suffices to consider those $\sig$ fixing
$\om_1$}. Theorem 1 is proved if we can show such a $\sig$ must equal the
identity.

A convenient way to exploit eq.(3a) is to use the {\it fusion rules}, which we
may take to be defined by {\it Verlinde's formula} [26]:
$$N_{\la,\mu}^{\nu}\Eqde\sum_{\be\in P_+^{r,k}}S_{\la,\be}\,{S_{\mu,\be}
\over S_{k\La_0,\be}}\,S^{*}_{\nu,\be}.\eqno(8a)$$
These $N_{\la,\mu}^\nu$ have a well-known geometric interpretation, but it
is irrelevant for our purposes.
Kac ([18]; see also p.288 of [17]) and Walton [28] used eq.(3a) to
reduce eq.(8a) to an expression involving finite-dimensional tensor
product decompositions. In particular, for any $\omega$ in the affine Weyl
group $W$ of $\A_r$, let $\epsilon(\omega)$ denote the sign of $\omega$,
and $\omega. \lambda=\omega(\la+\rho)-\rho$. Let $mult_{\bar{\la}
\otimes \bar{\mu}}(\bar{\nu})$ denote the multiplicity (i.e.\ the
Littlewood-Richardson coefficient) of the $A_r$-module
$\bar{L}(\bar{\nu})$ in $\bar{L}(\bar{\la})\otimes \bar{L}(\bar{\mu})$. Then
we find
$$N_{\la,\mu}^{\nu}=\sum_{\omega\in W}
\eps(\omega)\,\,mult_{\bar{\la}\otimes \bar{\mu}}\, (\overline{\omega.{\nu}}
).\eqno(8b)$$

We will use in Section 3 the following facts, obvious from eqs.(8a),(4d),(4e):
$$\eqalignno{N_{\la,k\La_0}^\nu=&\,\delta_{\la,\nu},\quad \forall \la,\nu\in
P_+^{r,k};&(8c)\cr
N_{C_r^cJ_r^a\la,C_r^cJ_r^b\mu}^{C_r^cJ_r^{a+b}\nu}=&\,N_{\la,\mu}^\nu,\quad
\forall\la,\mu,\nu\in P_+^{r,k},\ \forall a,b,c\in{\bf Z}.&(8d)\cr}$$

{}From eqs.(5a),(8a), $\sig$ must be
a symmetry of the fusion rules:
$$N_{\la,\mu}^\nu=N_{\sig\la,\sig\mu}^{\sig \nu}.\eqno(9a)$$
It is thus natural to look at the fusions involving $\om_1$: from eq.(8b) we
find
$$N_{\la,\om_1}^\nu=\left\{\matrix{1&{\rm if}\ \nu\in\la+\{\La_1-
\La_0,\La_2-\La_1,\ldots,\La_0-\La_{r}\}\cr 0&{\rm otherwise}\cr}
\right..\eqno(9b)$$

Now suppose $\sig\om_j=\om_j$ for all $1\le j\le i<r$. From eq.(9b), $N_{\om_i,
\om_1}^\nu=1$ iff either $\nu=\om_{i+1}$ or (for $k>1$) $\nu=\om_i'\Eqde
(k-2)\La_0+\La_1+\La_i$. This means, from eq.(9a), that $\sig\om_{i+1}$
must equal either $\om_{i+1}$ or $\om'_i$. But
$$(\om_{i}'+\rho\,|\,\om_{i}'+\rho)-(\om_{i+1}+\rho\,|\,\om_{i+1}+\rho) =
2i+2\not\equiv 0\quad ({\rm mod}\ 2\bar{k}).$$
Therefore by eq.(5b), $\sig\om_{i+1}=\om_{i+1}$.

Thus by induction $\sig$ must fix each $\om_1,\ldots,\om_r$. From this
we can now complete the proof of the theorem, with the following
observation:

\medskip\noindent{\bf Proposition 3}\quad {\it Suppose we have weights $\la,
\la'\in P_{+}^{r,k}$ such that
$${S_{\om_i,\la}\over S_{k\La_0,\la}}={S_{\om_i,\la'}\over S_{k\La_0,\la'}}
\eqno(10a)$$
for all $i$. Then} $\la=\la'$.

\noindent{\it Proof.}\quad We know (see Ch.VI, \S 3.4, Th.1 of [3]) that the
Weyl character $ch_{\overline{\be}}$ of the $\bar{g}$-module $\bar{L}(
\overline{\be})$ of any finite dimensional Lie algebra $\bar{g}$ can
be written as a polynomial $p_{\bar{\be}}$ of $ch_{\overline{\om_1}},\ldots,
ch_{\overline{\om_{r}}}$. From eq.(3a) we thus get
$${S_{\be,\ga}\over S_{k\La_0,\ga}}=p_{\bar{\be}}({S_{\om_1,\ga}\over
S_{k\La_0,\ga}},
\ldots,{S_{\om_{r},\ga}\over S_{k\La_0,\ga}}),\qquad\forall \beta,\gamma\in
P_+^{r,k}.\eqno(10b)$$
Equations (10a),(10b) together tell us that in fact
$${S_{\be,\la}\over S_{k\La_0,\la}}={S_{\be,\la'}\over S_{k\La_0,\la'}}
\qquad\forall\be\in P_+^{r,k}.$$
Multiplying it by $S^*_{\beta,\la}$ and summing over all $\beta$ forces
$\la=\la'$ by unitarity of $S$. \qquad QED\medskip

Choose any $\la\in P_{+}^{r,k}$ and put $\la'=\sig\la$. Then from
eq.(5a) and
using $\sig(k\La_0)=k\La_0$ and $\sig\om_i=\om_i$ we find that eq.(10a)
is satisfied. Then Prop.\ 3 tells us $\la=\la'$. In other words, $\sig$
must be the identity, and Theorem 1 is proved.

\bigskip
\noindent{\bf 3 \quad The automorphism invariants of $(A_{r_1}\oplus\cdots
\oplus A_{r_s})^{(1)}$} \medskip

The affinization  of reductive Lie algebras is discussed e.g.\ in \S\S
12.9-12.10 of [17]. Of course most quantities for semi-simple Lie algebras
can be built up in a straightforward way from those for the
simple ones. But this is not true for automorphism
invariants, as we shall see. Knowing the list of automorphism
invariants for the affinization of simple algebras $\bar{g}_i$ helps very
little in
their classification for the affinization of $\bar{g}_1\oplus\cdots\oplus
\bar{g}_s$. Nevertheless, with some additional
complications the techniques developed in the last section can be applied to
the classification of automorphism invariants for $\A_{\vec{r}}$.
Only one case with $s>1$ was known previously: all $r_i=1$ [11].

Theorem 2 gives the classification and is the main result of this section.
We would have liked to find an explicit set of generators for the group
of automorphism invariants, but this seems to be more work than it is worth,
for general $r_i$. Instead we will limit ourselves here to two special
cases, which together form Theorem 3.

The level here is an $s$-tuple $k=(k_1,\ldots,k_s)$, each $k_i$ a positive
integer.
Write $r=\vec{r}=(r_1,\ldots,r_s)$. As before call $\r_i=r_i+1$,
$\bar{k}_i=k_i+\r_i$. The set of highest weights is
$$P_+^{r,k}\Eqde P_+(\A_{\vec{r}},k)\ {\rm mod}\ {\bf C}\delta
=\bigl\{\la=\sum_{i=1}^s\la_{(i)}\,|\,\la_{(i)}=\sum_{j=0}^{r_i}
\la_{(i)j}\,\La_j^i\in P_+^{r_i,k_i}\bigr\},$$
using obvious notation. The modular matrix $S$ will be
$$S^{r,k}_{\la,\mu}=\prod_{i=1}^sS^{r_i,k_i}_{\la_{(i)},\mu_{(i)}},$$
similarly for $T^{r,k}$. We will usually drop the superscripts `$r,k$'.
For $a=(a_1,\ldots,a_s)$, write
$$\eqalignno{J_r^a\la\Eqde &\sum_{i=1}^sJ_{r_i}^{a_i}\la_{(i)},&\cr
C_r^a\la\Eqde&\sum_{i=1}^sC_{r_i}^{a_i}\la_{(i)},&(11a)\cr}$$
and $C_r=C_{r_1}\cdots C_{r_s}$. Let $k\La_0=\sum_jk_j\La_0^j$, and
for each $1\le i\le s$, $1\le \ell\le r_i$, define
$$\om^i_\ell\Eqde \sum_{j\ne i}k_j\La_0^j+(k_i-1)\,\La_0^i+\La_\ell^i\in
P_+^{r,k}.$$
Finally, write $[\la]$ for the orbit $\{C_r^{a}J_r^{b}\la\,|\,a_i,b_i\in
{\bf Z}\}$, and $Q(\la)=S_{\la,k\La_0}/S_{k\La_0,k\La_0}$. Prop.\ 1 implies
$$Q(\la)=1\ {\rm iff}\ \la\in[k\La_0].\eqno(11b)$$

Examples of automorphism invariants are the $C_r^a$ of eq.(11a). Another
example is induced
by any permutation $\pi$ of the indices $\{1,\ldots,s\}$ with the property
that
$$k_i=k_{\pi i},\qquad r_i=r_{\pi i},\qquad\forall i;\eqno(12a)$$
the corresponding automorphism invariant is the map $\sig_\pi$ defined by
$$\sig_\pi\la=\sum_{i=1}^s\la_{(\pi i)}.\eqno(12b)$$

Find any integers $a_{ij}$, for all $1\le i,j\le s$, satisfying:
$$\eqalignno{a_{ij}\,\r_i/\r_j\in &\,{\bf Z}\qquad \forall i,j;&(13a)\cr
{2a_{ii}\over \r_i}+\sum_{j=1}^sk_j{a^2_{ij}\over \r_j}\equiv &\,\sum_{j=1}^s
k_j\,a_{ij}\quad ({\rm mod\ 2}),\qquad \forall i;&(13b)\cr
{a_{ij}\over \r_j}+{a_{ji}\over \r_i}+\sum_{\ell=1}^sk_\ell{a_{i\ell}\,
a_{j\ell}\over \r_\ell}
\equiv &\,0\ \quad({\rm mod}\ 1),\qquad \forall i,j.&(13c)\cr}$$
To any such matrix $a$, define a function $\sig_a$ on $P_+^{r,k}$ by
$$\sig_a\la=J_r^{t(\la)\, a}\la\Eqde \sum_{j=1}^s\prod_{i=1}^s
J_{r_j}^{ a_{ij}t(\la_{(i)})}\,\la_{(j)}.\eqno(13d)$$
Because of eq.(13a), $\sig_a$ will be a permutation, with inverse $\sig_{b}$
where $b_{ij}=\r_j a_{ji}/\r_i$.
Because of eqs.(13b),(13c), $\sig_a$ will satisfy eqs.(5) and be an
automorphism invariant. Note that
$$\pi_b\circ\pi_a=\pi_c,\quad{\rm where}\ c_{ij}\equiv a_{ij}+b_{ij}+
\sum_{\ell=1}^sk_\ell\,a_{i\ell}\,b_{\ell j}\quad ({\rm mod}\ \r_j).
\eqno(13e)$$

\noindent{{\bf Theorem 2}}\quad {\it $\sig$ is an automorphism
invariant of $\A_{\vec{r}}$ at level $k=(k_1,\ldots,
k_s)$ iff $\sig=\sig_\pi\circ C_r^{c}\circ
\sig_a$, where $\sig_\pi$ is defined in eq.(12b), $C_r^c$ in eq.(11a),
and $\sig_a$ in eq.(13d).}\medskip

\noindent{{\it Remark.}}\quad
Note that in general we have here neither $\sig^2=id.$ nor $\sig\circ\sig'=
\sig'\circ\sig$. Also, we have here many more solutions than could be built
by tensoring together $s=1$ automorphism invariants. Provided we demand
$c_i=c_i'=0$ when $r_i=1$ or $k_i\le 2$, and $\pi i=\pi' i=i$ when $k_i=1$,
then
$$\sig_\pi\circ C_r^c\circ\sig_a=\sig_{\pi'}\circ C_r^{c'}\circ \sig_{a'}
\Leftrightarrow\pi i=\pi' i,\ c_i\equiv c'_i\ ({\rm mod}\ 2),\ a_{ij}
\equiv a'_{ij}\ ({\rm mod}\ \r_j)\ \forall i,j.$$

\noindent{\it Proof.}\quad Write $Q_i=Q(\om^i_1)$.  For convenience reorder
the indices so that $Q_1\le Q_2\le\cdots\le Q_s$. We will begin by finding
an index $\pi i$ for each $i$, so that $\sig\om^i_1
\in[\om^{\pi i}_1]$. As in eq.(8a) put
$$N_{\la,\mu}^\nu\Eqde\sum_{\beta\in P^{r,k}_+}S^{r,k}_{\la,\beta}\,{
S^{r,k}_{\mu,\beta}\over S^{r,k}_{k\La_0,\beta}}\,S_{\nu,\beta}^{r,k*}=
\prod_{i=1}^sN_{\la_{(i)},\mu_{(i)}}^{\nu_{(i)}}.\eqno(14a)$$
Let $o_i(\la)$ equal the number of $0\le j\le r_i$ such
that $\la_{(i)j}>0$. Then from eqs.(9b),(8c),
$$\sum_{\nu\in P_+^{r,k}}N_{\la,\om^i_1}^\nu=o_i(\la),\qquad
\forall i\in\{1,\ldots,s\},\ \la\in P_+^{r,k}.\eqno(14b)$$

We will construct $\pi$ by induction on $i$. Suppose for all $j<i$, we
have a $\pi j$ such that $\sig\om^j_1\in[\om^{\pi j}_1]$.
By eqs.(8d),(14a),(14b),(5a) we know for these $j<i$ that
$$o_j(\la)=o_{\pi j}(\sig\la),\qquad \forall \la\in P_+^{r,k}.\eqno(14c)$$

If $Q_i=1$ (i.e.\ $k_i=1$), then both $\om^i_1,\sig\om^i_1\in[k\La_0]$ by
eqs.(5a),(11b), so put $\pi i=i$ there ($k_i=1$ is a special
-- albeit trivial -- case here because for those $i$, $\pi i$ is not fixed by
the  constraint $\sig\om^i_1\in[\om^{\pi i}_1]$). When $Q_i>1$,
$(\sig\om^i_1)_{(j)}\in [k\La_0]_{(j)}$ for all $j$ with $Q_j<Q_i$,
 by  eq.(14c). So by Prop.\ 1, we
must have $\sig\om^i_1\in[\om^\ell_1]$ for some $\ell$ with $Q_i=
Q_\ell$. Put $\pi i=\ell$.

This inductively defines $\pi i$ for all $i$, in such a way that
$\sig\om^i_1\in[\om^{\pi i}_1]$.
Next, let us show $\pi$ is a bijection. It suffices to show that $\pi$ is
onto. Choose any $1\le \ell\le s$ -- we may suppose $k_\ell>1$. Because
$\sigma^{-1}$ is also an automorphism invariant, there exists an $i$ such
that $\sig (C_r)^aJ_r^b\om^i_1=\om^\ell_1$ for some $a\in\{0,1\},b=(b_1,
\ldots,b_s)$.

Now $\sig$ commutes with $S$, so it also does with $C^{r,k}=(S^{r,k})^2$.
Also, by eqs.(5a),(11b) $\sig J_r^b(k\La_0)=J_r^c(k
\La_0)$ for some $c$. But because of eqs.(4d),(3c),(5a) we know
$$S_{J_r^b\la,\mu}=S_{\la,\mu}\,{S_{J_r^b(k\La_0),\mu}\over S_{k\La_0,\mu}}=
S_{\sig\la,\sig\mu}\,{S_{J_r^c(k\La_0),\sig\mu}\over S_{k\La_0 ,\sig\mu}}
=S_{J_r^c\sig\la,\sig\mu},\quad \forall \la,\mu\in P_+^{r,k}.$$
Therefore by unitarity of $S$, $\sig J^b_r=J^c_r\sig$. Hence $\sig\om^i_1=
J^{-c}_r(C_r)^a\om^\ell_1\in [\om^\ell_1]$, so $\pi i=\ell$, and
$\pi $ is a bijection.

Finally we want to show $\pi$ satisfies eq.(12a). Write $\sig\om^i_1=(C_r)^{
c_i} J_r^d\om^{\pi i}_1$, for some $c_i\in\{0,1\}$, $d=(d_1,\ldots,d_s)$. The
expression
$$1=N_{\om^i_1,\om^i_1}^\la=
N_{(C_r)^{c_i}J_r^d\om^{\pi i}_1,(C_r)^{c_i}J_r^d\om^{\pi i}_1}^{\sig\la}$$
requires from eqs.(14a),(8c),(8d),(9b) that $\la=\om^i_2$ or $\la=
\om^i_1{}'\Eqde \sum_{j\ne i}k_j\,\La_0^j+(k_i-2)\La_0^i+2\La_1^i$, and
that $\sig\la=(C_r)^{c_i}J_r^{2d}\om^{\pi i}_2$ or $\sig\la=(C_r)^{c_i}J^{2d}_r
\om^{\pi i}_1{}'$. Thus either:

\item{(i)} $\sig\om^i_2=(C_r)^{c_i}J_r^{2d}\om^{\pi i}_2$ and
$\sig\om^i_1{}'=(C_r)^{c_i}J_r^{2d}\om^{\pi i}_1{}'$; or

\item{(ii)} $\sig\om^i_2=(C_r)^{c_i}J_r^{2d}\om^{\pi i}_1{}'$ and
$\sig\om^i_1{}'=(C_r)^{c_i}J_r^{2d}\om^{\pi i}_2$.

But from eqs.(5b),(4c) and using $\sig\om^i_1=(C_r)^{c_i}J_r^d\om^{\pi i}_1$,
case (i) is seen to require that
$${\r_i\pm 1\over \overline{k}_i}\equiv {\r_{\pi i}\pm 1\over
\overline{k}_{\pi i}}\quad ({\rm mod}\ 1);\eqno(15a)$$
while case (ii) requires
$${\r_i\pm 1\over \overline{k}_i}\equiv {\r_{\pi i}\mp 1\over
\overline{k}_{\pi i}}\quad ({\rm mod}\ 1).\eqno(15b)$$
Since $2\le \r_j\le \bar{k}_j-1$ for all $j$, eq.(15a) forces eq.(12a) while
 eq.(15b) has no solution.

Thus $\pi$ satisfies eq.(12a). Replacing $\sig$ by $\sig_{\pi^{-1}}\circ
\sig$, we may assume $\sig\om^i_1\in[\om^i_1]$ for all $i$. Likewise,
by replacing $\sig$ with $C_r^c\sig$ for some $c=(c_1,\ldots,c_s)$ (we can
require $c_i=0$ when $r_i=1$ or $k_i\le 2$), we may
suppose, for each $i$, that
$$\sig\om^i_1=J_r^{a_i}\om^i_1,\qquad {\rm where}\ a_i=(a_{i1},\ldots,
a_{is}).\eqno(16a)$$
Eqs.(5a),(4d) and (5b),(4c) give us eqs.(13c),(13b) respectively. Using a
calculation given earlier this proof,
$\sig^{-1}$ will be an automorphism invariant satisfying eq.(16a) for
some matrix $b$ in place of $a$. Looking at $S_{\om^i_1,\sig^{-1}
\om^j_1}=S_{\sig\om^i_1,\om^j_1}$ we find from eq.(4d)
$$b_{ji}/\r_i\equiv a_{ij}/\r_j\ \quad ({\rm mod}\ 1),\eqno(16b)$$
thus eq.(13a) will also be satisfied. This means $\sig_a$ defines an
automorphism
invariant; replacing $\sig$ with $\sig_a^{-1}\circ\sig$, we may then assume
all $a_{ij}=0$. It suffices now to prove any such $\sig$ must be the identity.

By exactly the same induction argument used in Section 2, together with
eqs.(8c),(14a), (9b) the expression
$$1=N_{\om^j_i,\om^j_1}^{\om^j_{i+1}}=N_{\om^j_i,\om^j_1}^{\sig\om^j_{i+1}}$$
gives only two possibilities for $\sig\om^j_{i+1}$. As before, only
$\sig\om^j_{i+1}=\om^j_{i+1}$ satisfies eq.(5b). This argument shows that
$\sig$ must fix all $\om^\ell_m$.
The proof of Prop.\ 3 now carries over without change, and we find that
indeed $\sig\la=\la$ for all $\la\in P_+^{r,k}$.  \qquad QED\medskip

As is done in Section 2 for the special case $s=1$, it should be possible
to find a complete set of generators for the group of automorphism invariants
of $\A_{\vec{r}}$ at fixed level $k$, as well as compute its order.
But both of these will be messy, depending on $r$ and $k$ in
a more complicated way than is the case for $s=1$. We will limit
ourselves here to two simple observations (see Theorem 3 below).

But first, define four new types
of automorphism invariants $\sig[J_r^m]$, $\sig[p;\ell,m]$, $\sig[\ell,m,n]$
and $\sig[\ell,m,n,o]$ as follows.

\item{(i)} Choose any $m=(m_1,\ldots,m_s)$. Define $\tilde{k}_i
=\bar{k}_i$ or $k_i$ as in eq.(6a). Put $u=\sum_i\tilde{k}_im_i\,(\r_i-m_i)/
\r_i$, $\bar{m}_i=\r_i/gcd\{m_i,\r_i\}$, and $N=lcm_i\{
\bar{m}_i\}$. Suppose $Nu$ is even, and $Nu/2$ is coprime to $N$. Then
there exists an integer $v$ such that $vNu/2\equiv 1$ (mod $N$).
Define $\sig[J_r^m]=\sig_a$, where $a_{ij}=vNm_jm_i/\r_i$ (this is the
immediate generalization of eq.(6b)).

\item{(ii)} Choose any prime $p$, and any
indices $1\le \ell<m \le s$ such that $p$ divides $\r_\ell$, $k_\ell$,
$\r_m$ and $k_m$. If $p=2$ we need the additional constraint that 8 divides
both $\r_\ell k_\ell$ and $\r_mk_m$.
Define $\sig[p;\ell,m]$ to be $\sig_a$, where $a_{ij}=0$
for all $i,j$ except for $a_{\ell m}=\r_m/p$, $a_{m \ell}=-\r_\ell/p$.

\item{(iii)} Find
distinct indices $\ell,m,n$ with $k_\ell,k_m$ both odd, $k_n\equiv 0$
(mod 4), and $\r_\ell\equiv \r_m\equiv\r_n\equiv 2$ (mod 4), such that
$k_\ell\r_\ell+k_m\r_m\equiv 0$ (mod 8). Define $\sig[\ell,m,n]$ to be
$\sig_a$, where $a_{ij}=0$ except for $a_{\ell n}=a_{mn}=\r_n/2$, $a_{n\ell}
=\r_\ell/2$, and $a_{nm}=\r_m/2$.

\item{(iv)} Find distinct indices $\ell,m,n,o$, with
$k_\ell,k_m,k_n,k_o$ all odd, and $\r_\ell\equiv\r_m\equiv
\r_n\equiv\r_o\equiv 2$ (mod 4), such that $k_m\r_m+k_n\r_n\equiv k_\ell
\r_\ell+k_o\r_o\equiv 0$ (mod 8).
Define $\sig[\ell,m,n,o]$ to be $\sig_a$, where $a_{ij}
=0$ except for $a_{\ell m}=a_{om}=\r_m/2$, $a_{\ell n}=a_{on}=\r_n/2$,
$a_{m\ell}=a_{n\ell}=\r_\ell/2$ and $a_{mo}=a_{no}=\r_o/2$.

\medskip\noindent{{\bf Theorem 3}}\quad {\it (a) If
$gcd\{2k_i,\r_i\}=1$ for each $i$, then the $\sig[J_r^m]$ of (i) above,
together with $C_r^c$ and $\sig_\pi$, generate
the group of all automorphism invariants of $\A_{\vec{r}}$.}

\item{{\it (b)}} {\it If each $\r_i$ is square-free, then the group of
automorphism invariants of $\A_{\vec{r}}$ is generated by
$C_r^c$, $\sig_\pi$, and the automorphism invariants $\sig[J_r^m]$,
$\sig[p;\ell,m]$, $\sig[\ell,m,n]$, and $\sig[\ell,m,n,o]$ defined in
(i)-(iv) above.}\medskip

To prove Theorem 3, note that it suffices to fix a prime $p$ and to
consider all solutions $a$ to eqs.(13a)-(13c) when the $\r_i$
are all powers of $p$. Both Theorems 3(a),3(b) follow
inductively from the following observation. Suppose $\exists \ell$ such that
$${\r_\ell\over gcd\{\r_\ell,a_{\ell\ell}\}}\ge {\r_j\over gcd\{\r_j,a_{j\ell}
\r_j/\r_\ell\}},\qquad \forall j\ {\rm with}\ a_{j\ell}\not\equiv 0\ ({\rm
mod}\ \r_\ell).\eqno(17)$$
Then we may take $m_j=a_{j\ell}\r_j/\r_\ell$: $\sig[J_r^m]$ will be an
automorphism invariant, and $\sig_b\Eqde \sig_a\circ(\sig[J_r^m])$ will have
$b_{\ell j}=b_{j\ell}=0$ $\forall j$, by eq.(13e). This means
$(\sig_b\la)_{(\ell)}=\la_{(\ell)}$ $\forall \la\in P_+^{r,k}$.

To prove Theorem 3(a), let $\ell$ satisfy $\r_\ell\ge \r_j$ $\forall j$ with
$a_{j\ell}\not\equiv 0$ (mod $\r_\ell$). If
$gcd\{\r_\ell,a_{\ell\ell}\}\ne 1$, then replace $\sig_a$ with $\sig_{a'}=
\sig_a \circ(\sig[J_r^m])$ where $m_{i}=\delta_{i\ell}$ -- then eq.(13e)
tells us $gcd\{\r_\ell,a'_{\ell\ell}\}=1$, so eq.(17) is satisfied.

Theorem 3(b) follows from similar arguments. The difficult case there is $p=2$,
which is worked out explicitly in [11].

\bigskip\noindent{\bf 4\quad Conclusion}\medskip

We begin the paper by reviewing the problem of classifying conformal field
theories -- in particular what are called their partition functions. An
important class
of these are the {\it automorphism invariants} defined in eqs.(5). This paper
classifies all automorphism invariants corresponding to the affinization of
$A_{r_1}\oplus\cdots\oplus A_{r_s}$. This is a necessary and
major step toward the full classification of all conformal field theories
corresponding to those affine algebras. This result, especially the
simplicity of its proof, is a strong hint that the full classification
should be possible.

We find here that all such automorphism invariants are related in a fairly
simple way to the symmetries of the corresponding extended Dynkin diagrams
(this will not always be true for other affine algebras [27,12]). For $\A_r$,
they all commute and have order 2 (this is not true for most other $g$ --
e.g.\ $(A_{r_1}\oplus\cdots\oplus A_{r_s})^{(1)}$ for $s\ge 2$).

A simple sketch of the proof is given at the beginning of Section 2.
Our arguments are mostly combinatorial, though a number of
algebraic results are required. A remarkable fact used
repeatedly is eq.(3a), which says that the ratio of certain entries
of the modular matrix $S$ of $g$ equals a certain value of a Weyl character
of $\bar{g}$. One place this is exploited is
in the rewriting of Verlinde's formula in terms of finite-dimensional
tensor product multiplicities.

The arguments here should extend to all other affine algebras [12],
but $A_r$ is notoriously well-behaved for a Lie algebra so
there will be some additional complications,
particularly at small levels. The arguments should also lift
to more general conformal field theory classifications (namely,
what are called the ${\cal D}$-type invariants and ${\cal E}_7$-type
exceptionals [11]).

Considering the diverse applications of the representation theory of Kac-Moody
algebras, it can perhaps be hoped that the work in this paper -- namely
the classification of all symmetries of the Kac-Peterson modular matrices --
will also find application outside the scope of conformal field theory.

\bigskip\noindent{{\it Acknowledgements.}} \quad I am especially grateful to
Philippe Ruelle
for convincing me that this problem is solvable,
and to Mark Walton for his patient explanations of fusion rules.
I also thank Antoine Coste and Jean-Bernard Zuber for general
information, and Ron Wang for many helpful stylistic remarks. The hospitality
 of IHES is also appreciated.\bigskip

\noindent{{\bf References}}\medskip

\item{{\bf 1.}} Altschuler, D., Lacki, J., Zaugg, Ph.: The affine Weyl group
and modular invariant partition functions. Phys.\ Lett.\ {\bf B205}, 281-284
(1988)

\item{{\bf 2.}} Bernard, D.: String characters from Kac-Moody automorphisms.
Nucl.\ Phys.\ {\bf B288}, 628-648 (1987)

\item{{\bf 3.}} Bourbaki, N.: Groupes et Alg\` ebres de Lie, Chapitre IV-VI.
Paris: Hermann 1968

\item{{\bf 4.}} Bouwknegt, P., Nahm, W.: Realizations of the exceptional
modular
invariant $A^{(1)}_1$ partition functions. Phys.\ Lett.\ {\bf B184}, 359-362
(1987)

\item{{\bf 5.}} Cappelli, A., Itzykson, C., Zuber, J.-B.: The A-D-E
classification
of $A^{(1)}_1$ and minimal conformal field theories. Commun.\ Math.\
Phys.\ {\bf 113}, 1-26 (1987)

\item{{\bf 6.}} Coste, A., Gannon, T.: Remarks on Galois symmetry in rational
conformal field theory. Phys.\ Lett.\ {\bf B323}, 316-321 (1994)

\item{{\bf 7.}} Coste, A., Itzykson, C., Lascoux, J., Ruelle, Ph.
(private communications)

\item{{\bf 8.}} Fuchs, J.: Simple WZW currents. Commun.\ Math.\ Phys.\ {\bf
136}, 345-356 (1991)

\item{{\bf 9.}} Gannon, T.: WZW commutants, lattices, and level-one partition
functions. Nucl.\ Phys.\ {\bf B396}, 708-736 (1993)

\item{{\bf 10.}} Gannon, T.: The classification of affine SU(3) modular
invariant partition functions. Commun.\ Math.\ Phys.\ {\bf 161}, 233-264
(1994); The classification of SU(3) modular invariants revisited. IHES
preprint P/94/32 (hep-th/9404185)

\item{{\bf 11.}} Gannon, T.: Towards a classification of SU(2)$\oplus\cdots
\oplus$SU(2) modular invariant partition functions. J.\ Math.\
Phys.\ (to appear) (hep-th/9402074)

\item{{\bf 12.}} Gannon, T., Ruelle, Ph., Walton, M.\ A. (in preparation)

\item{{\bf 13.}} Gato-Rivera, B., Schellekens, A.\ N.: Complete classification
of simple current automorphisms. Nucl.\ Phys.\  {\bf B353}, 519-537 (1991)

\item{{\bf 14.}} Gepner, D., Witten, E.: String theory on group manifolds.
Nucl.\ Phys.\ {\bf B278}, 493-549 (1986)

\item{{\bf 15.}} Goddard, P., Kent, A., Olive, D.: Unitary representations
of the Virasoro and super-Virasoro algebras. Commun.\ Math.\ Phys.\
{\bf 103}, 105-119 (1986)

\item{{\bf 16.}} Itzykson, C.: Level one Kac-Moody characters and modular
invariance. Nucl.\ Phys.\ (Proc.\ Suppl.) {\bf B5}, 150-165 (1988)

\item{{\bf 17.}} Kac, V.\ G.: Infinite Dimensional Lie Algebras,
3rd ed. Cambridge: Cambridge University Press 1990

\item{{\bf 18.}} Kac, V.\ G.: talk at Canadian Mathematical Society Meeting on
Lie algebras and Lie groups (CRM, Universit\'e de Montr\'eal, August 1989),
unpublished

\item{{\bf 19.}} Kac, V.\ G., Peterson, D.: Infinite-dimensional Lie algebras,
theta functions and modular forms. Adv.\ in Math.\ {\bf 53}, 125-264 (1984)

\item{{\bf 20.}} Kac, V.\ G., Wakimoto, M.: Modular and conformal constraints
in respresentation theory of affine algebras. Adv.\ in Math.\ {\bf 70},
156-236 (1988)

\item{{\bf 21.}} Kass, S., Moody, R.\ V., Patera, J., Slansky, R.: Affine Lie
Algebras, Weight Multiplicities, and Branching Rules, Vol.\ 1. Berkeley,
Los Angeles, Oxford: University of California Press 1990

\item{{\bf 22.}} Koblitz, N., Rohrlich, D.: Simple factors in the Jacobian of
a Fermat curve. Can.\ J.\ Math.\ {\bf XXX}, 1183-1205 (1978)

\item{{\bf 23.}} Ruelle, Ph., Thiran, E., Weyers, J.: Implications of an
arithmetic symmetry of the commutant for modular invariants. Nucl.\ Phys.\
{\bf B402}, 693-708 (1993)

\item{{\bf 24.}} Ruelle, Ph.: Automorphisms of the affine SU(3) fusion rules.
Commun.\ Math.\ Phys.\ {\bf 160}, 475-492 (1994)

\item{{\bf 25.}} Schellekens, A.\ N., Yankielowicz, S.: Modular invariants
from
simple currents. An explicit proof. Phys.\ Lett.\ {\bf B227}, 387-391 (1989)

\item{{\bf 26.}} Verlinde, E.: Fusion rules and modular transformations in 2D
conformal field theory. Nucl.\ Phys.\ {\bf B300 [FS22]}, 360-376 (1988)

\item{{\bf 27.}} Verstegen, D.: New exceptional modular invariant partition
functions for simple Kac-Moody algebras. Nucl.\ Phys.\ {\bf B346}, 349-386
(1990)

\item{{\bf 28.}} Walton, M.\ A.: Algorithm for WZW fusion rules:
a proof. Phys.\ Lett.\ {\bf B241}, 365-368 (1990)

\item{{\bf 29.}} Witten, E.: Non-abelian bosonization in two dimensions.
Commun.\ Math.\ Phys.\ {\bf 92}, 455-472 (1984)

\item{{\bf 30.}} J.-B.\ Zuber, talk given at ICMP, Paris, 1994

\end